\documentclass[journal]{IEEEtran}
\usepackage{amsmath}
\usepackage{amssymb}
\usepackage[dvips]{graphicx}
\usepackage{color}
\usepackage{subfigure}
\usepackage{bm}

\newtheorem{RK}{Remark}

\newtheorem{define}{Definition}

\DeclareMathOperator*{\argmax}{argmax}
\DeclareMathOperator*{\diag}{diag}
\newcommand{\nn}{\nonumber}
\newcommand{\mc}{\mathcal}
\newcommand{\mb}{\mathbf}
\newcommand{\ms}{\mathsf}

\newcommand{\mbb}{\mathbb}
\newcommand{\bs}{\boldsymbol}
\DeclareMathOperator{\rank}{rank}

\ifCLASSINFOpdf
\else
\fi
\begin{document}
%
\title{On PMU Location Selection for Line Outage Detection in Wide-area Transmission Networks}
%
%
%

\author{Yue~Zhao,~\IEEEmembership{Member,~IEEE,}
        Andrea~Goldsmith,~\IEEEmembership{Fellow,~IEEE,}
        and~H.~Vincent~Poor,~\IEEEmembership{Fellow,~IEEE}

\thanks{This work was supported in part by the NSF under Grant CNS-09-05086, and in part by the DTRA under Grant HDTRA1-08-1-0010.}%
\thanks{Y. Zhao is with the Department of Electrical Engineering, Princeton University, Princeton, NJ, 08544 USA, and with the Department of Electrical Engineering, Stanford University, Stanford, CA, 94305 USA (e-mail: yuez@princeton.edu).}%
\thanks{A. Goldsmith is with the Department of Electrical Engineering, Stanford University, Stanford, CA 94305 USA (e-mail: andrea@stanford.edu)}%
\thanks{H. V. Poor is with the Department of Electrical Engineering, Princeton University, Princeton, NJ 08544 USA (e-mail: poor@princeton.edu).}}%
\maketitle

\begin{abstract}
The optimal PMU locations to collect voltage phase angle measurements for detecting line outages in wide-area transmission networks are investigated. The problem is established as one of maximizing the minimum distance among the voltage phase angle signatures of the outages, 
which can be equivalently formulated as an integer programming problem. Based on a greedy heuristic and a linear programming relaxation, a branch and bound algorithm is proposed to find the globally optimal PMU locations. Using this algorithm, the optimal trade-off between the number of PMUs and the outage detection performance is characterized for IEEE 14, 24 and 30 bus systems. The algorithm is shown to find the globally optimal PMU locations in a small number of iterations. It is observed that it is sufficient to have roughly one third of the buses providing PMU measurements in order to achieve the same outage detection performance as with all the buses providing PMU measurements.

\end{abstract}


%
\IEEEpeerreviewmaketitle

\section{Introduction}
In high voltage transmission networks, lack of wide-area situational awareness has been one of the major causes of large-scale blackouts \cite{black04}. One of the reasons for the lack of situational awareness has been the limitations of the conventional sensors and SCADA (supervisory control and data acquisition) systems. Phasor measurement units (PMUs), as compared to conventional sensors, are able to provide GPS-synchronized, more accurate and temporally much denser measurements of voltage and current phasors \cite{PMUorig}. Thus, the deployment of PMUs for grid sensing is widely considered to be a major driving force for improving the reliability of transmission networks.

There has been considerable research investigating how PMU measurements can be exploited in various tasks for achieving a reliable grid, including outage detection \cite{lineoutage, doubleline, Hao}, state estimation \cite{AburSE05, obs93, obs04}, stability analysis \cite{DySec07}, etc. As the cost of installing and networking PMUs is relatively high, two of the major questions are i) what is the necessary or sufficient number of PMUs to use for achieving good performance in all these tasks, and ii) where are the optimal locations to collect PMU measurements? For achieving full network observability, it was estimated that about one third of the buses need to provide PMU measurements \cite{obs93}. However, it is unclear whether this estimate applies to \emph{other} kinds of tasks including outage detection. Moreover, for different tasks, the objectives for which the PMU measurements are used differ considerably, and the corresponding optimal PMU locations can vary greatly.

In this paper, we focus on the problem of using PMU measurements of \emph{voltage phase angles} for real time detection of \emph{line outages} over wide areas (including the ones that occur at places where no PMU measurements are available at the control center). We assume that network state estimation during normal conditions is available over relatively slow timescales, (e.g., using NERC System Data Exchange (SDX) \cite{SDX06}), which provides the pre-outage base case system parameters and states. In \cite{lineoutage, doubleline}, and \cite{Hao}, it was demonstrated that many line outages can be detected with only a \emph{subset} of the buses providing PMU measurements. However, a comprehensive understanding of the outage detection performance given arbitrary constraints on the number of PMUs to use is left open.

We address this open question of characterizing the optimal trade-off between the number of PMUs used and the outage detection performance. The central problem is finding the \emph{optimal locations} at which to collect PMU measurements, which is an NP hard combinatorial optimization. We consider \emph{both non-adaptive and adaptive} PMU location selections.

We first define a set of voltage phase angle signatures associated with the potential outages. Intuitively, the outage detection performance depends on the degrees of separation among the outage signatures. We employ the minimum distance among the signatures as an indicator of their degrees of separation. Accordingly, we introduce the max-min distance criterion for optimizing the PMU locations given any constraint on the number of PMUs. A simple greedy heuristic is developed that provides lower bounds on the global optimum. Next, we show that finding the optimal PMU locations can be equivalently formulated as an integer programming problem (IP), which allows a relaxation as a linear program (LP) that provides upper bounds on the global optimum. We then develop a branch and bound algorithm for PMU location selection: at each iteration, we strategically fix a bus to either provide or not provide PMU measurements, and compute new lower and upper bounds using the greedy heuristic and LP relaxation. 

To understand the optimal trade-off between the number of PMUs and the minimum distance among the outage signatures, we apply the proposed branch and bound algorithm on detecting all single line outages in IEEE 14, 24 and 30 bus systems \cite{MatP11}, and show that the globally optimal PMU locations can be found in a small number of iterations. The optimal trade-offs show that having roughly \emph{one third} of the buses to provide PMU measurements is sufficient for achieving the same outage detection performance as with all the buses providing PMU measurements.

\section{Problem Formulation}
We consider a power transmission network with $N$ buses and $L$ transmission lines. We denote the set of all the buses by $\mc{N}=\{1,2,\ldots,N\}$, and the set of all the lines by $\mc{L} = \{1,2,\ldots,L\}$. We study the problem of \emph{line outage detection}, where an outage event corresponds to the loss of a subset of lines in $\mc{L}$. 
Clearly, there are in total $2^L - 1$ different possible outage events, which correspond to all the non-empty subsets of $\mc{L}$. In practice, however, there are several reasons for us to consider a much smaller subset of the outage events:
\begin{itemize}
\item In case a single line outage results in the overloading and overheating of another line, it often takes minutes after the initial single line outage for the overheated line to trip. Thus, within a few tens of seconds, it is very unlikely to see many line outages happen together due to their probabilistic independence.
\item Not all the outage events have the same level of impact on the grid. For example, some line outages do not lead to the overheating of other lines, while others do. The latter have a greater impact on the grid as they may cause cascading failures.
\item As we would like to detect and locate outages in real time with good detection performance, the limitations of computational power may forbid us to consider too many outage events.
\end{itemize}
Based on the above considerations and the real world situation, a set of \emph{outage events of interest} shall be selected, which we denote by $\mc{E} = \{E_1, E_2, \ldots, E_K\}$. For example, in \cite{lineoutage, doubleline}, $\mc{E}$ contains all single and double line outages, respectively.

Given an outage event set $\mc{E}$, we augment it by the \emph{non-outage event} (i.e., normal condition) $E_0$, and denote the augmented set by $\bar{\mc{E}} = \{E_0,E_1,\ldots, E_K\}$. We investigate the following problem:
\begin{itemize}
\item The grid starts in normal condition $E_0$.
\item Either the grid stays in normal condition $E = E_0$, or an arbitrary outage event $E = E_k\in\mc{E}$ occurs.
\item We make a detection decision $\hat{E}\in\bar{\mc{E}}$. The detection is successful if $\hat{E} = E$, otherwise not.
\end{itemize}

\begin{figure}[tb]
  \centering
  \includegraphics[scale=0.45]{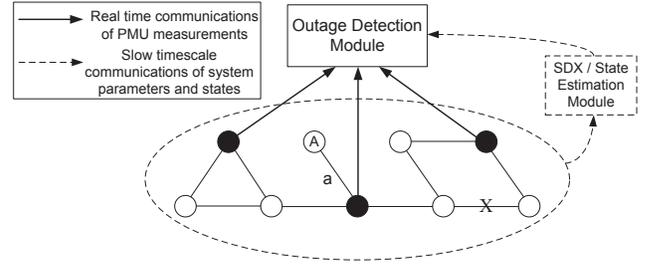}
  \caption{System diagram for the problem of line outage detection using PMU measurements of voltage phase angles. Buses and transmission lines are represented by circles and connecting segments. The solid circles are the buses where PMU measurements are communicated to the outage detection module in the control center. The solid links represent real time communications of the measurements. The dashed SDX / state estimation module and links provide the slow timescale updates of the system parameters and states to the outage detection module. The loss of line $X$ is an example of single line outage.}
  \label{sysdiag}
\end{figure}

We detect outages by observing the \emph{voltage phase angles at a subset of the buses}, where the phasor measurements are provided by PMUs. A system diagram for this problem is depicted in Figure \ref{sysdiag}. As outages shall be detected as quickly as possible to prevent cascading effects, we aim at detecting an outage within about ten seconds after its occurrence. 
We note that the phase angles typically stabilize within a few seconds after a line outage occurs \cite{lineoutage}, and we assume that the power injections of the network remain the same within these few seconds. This is a reasonable assumption as power injections change relatively slowly. We make use of the stabilized phase angles before and after an outage to make detection decisions.

\subsection{Pre-computing the Outage Signatures} \label{sigsec}
For each event $E_k\in\mc{E}$, we denote the corresponding stabilized phase angles at all the buses by $\bm{\theta}^{(k)}\in\mbb{R}^{N\times1}$. We employ the DC power flow model to compute $\bm{\theta}^{(k)}$ as follows:
\begin{align}
\bm{P} = \bm{B}^{(k)}\bm{\theta}^{(k)}, k=0,1,\ldots,K, \label{dcpf}
\end{align}
where $\bm{P}\in\mbb{R}^{N\times1}$ denotes the power injections, and 
$\bm{B}^{(k)}\in\mbb{R}^{N\times N}$ is determined by the line status information after event $E_k$ occurs \cite{WW96}. We note the following properties of $\mb{B}^{(k)}$:

\begin{itemize}
\item ~\vspace{-12pt}
\begin{equation}\bm{B}^{(k)}\bs{1} = \bs{0}, \label{LapRK}\end{equation}
\vspace{-18pt}

and we can arbitrarily add the same constant phase angle to every entry of $\bm{\theta}$, without having any change of $\bm{P}$.
\item $\rank(\bm{B}^{(k)}) = N - 1$ if and only if the grid is connected. In the case where the grid contains islands, $\rank(\bm{B}^{(k)}) = N - C$ where $C$ is the number of connected components (islands) of the grid.
\end{itemize}

We note that if an outage creates islands in the grid, the balance of the power injections may not be satisfied by the same $\bm{P}$, and \eqref{dcpf} may not hold anymore. For example, in Figure \ref{sysdiag}, the loss of line $a$ will create a single bus island $A$ which may not be self-balanced (i.e., the pre-outage power injection at $A$ may not be zero). We restrict our consideration to the cases where the post-outage grid remains connected, and leave the islanding cases as future work.

In practice, \emph{noisy} versions of the pre-outage and post-outage phase angle vectors are observed:
\begin{align}
\bm{\theta} = \bm{\theta}^{(0)} + \bm{z}^{0}, ~ \tilde{\bm{\theta}} = \bm{\theta}^{(\kappa)} + \bm{z}^1, \label{noisyv}
\end{align}
where $E_\kappa,\kappa\in{0,1,\ldots,K}$ is the actual outage (or non-outage if $\kappa = 0$) that occurred, and $\bm{z}^0,\bm{z}^1$ are observation noise vectors that account for the errors in the measurement data and the system parameter data. In the following sections, we assume that the noises at all the buses are i.i.d.. We will generalize this assumption in Section \ref{gennoise}. The task of line outage detection can then be formulated as the following hypothesis testing problem:

\emph{From observing the pre-outage and post-outage phase angle vectors $\bm{\theta}, \tilde{\bm{\theta}}$ \eqref{noisyv}, identify which event $E_k\in\bar{\mc{E}}$ has occurred}.

For this hypothesis testing problem, $\bm{\theta}^{(0)}, \bm{\theta}^{(1)},\ldots,\bm{\theta}^{(K)}$ can be viewed as the \emph{signatures} of the events $E_0,E_1,\ldots,E_K$. To identify any event $E_k\in\bar{\mc{E}}$, all the signatures $\{\bm{\theta}^{(k)}, k=0,1,\ldots,K\}$ shall be collected before an outage occurs.

From \eqref{dcpf}, we pre-compute the signatures by
\begin{align}
\bm{\theta}^{(k)} = {\bm{B}^{(k)}}^+\bm{P},k=0,1,\ldots,K, \label{precomp}
\end{align}
where ${\bm{B}^{(k)}}^+$ is the pseudoinverse of $\bm{B}^{(k)}$. 
We assume that the knowledge of the power injections $\bm{P}$ and the normal condition $\bm{B}^{(0)}$ matrix are available from system-wide data sources (e.g., the NERC SDX \cite{SDX06}, or other state estimation mechanisms that operate over relatively slow timescales, cf. Figure \ref{sysdiag}.) Note that, with the knowledge of the open line indices for event $E_k$, the post-outage $\bm{B}^{(k)}$ matrices $(k=1,2,\ldots,K)$ can be derived from the normal condition $\bm{B}^{(0)}$ matrix. As a result, while the grid is working under normal conditions, the signatures $\{\bm{\theta}^{(k)}, k=0,1,\ldots,K\}$ can be pre-computed by \eqref{precomp} in preparation for detecting the potential line outages in $\mc{E}$. Intuitively, \emph{the more separated the signatures are from each other, the better detection performance can be achieved.} In later sections, we characterize the separation by the minimum distance among the signature set. 

From \eqref{LapRK}, a phase angle vector $\bm{\theta}$ remains functionally equivalent by adding the same constant to all its entries. Without loss of generality (WLOG), we can choose any one of the $N$ buses as the \emph{reference bus}, 
and subtract the phase angle at the reference bus from the phase angles at all the $N$ buses, such that the phase angle at the reference bus is always kept at zero. In principle, it does not matter which bus among the $N$ buses is chosen to be the reference bus. However, it does matter in the case when only a \emph{subset} of the buses' phase angles are available. This is because once we choose a bus as the reference bus, we implicitly assume that the phase angle at this bus is available.

\subsection{PMU Location Selection} \label{pmuloc}
We now establish the main problem of this paper, namely, PMU location selection. In \eqref{noisyv}, synchronous and timely updates of the complete phase angle vector require \emph{all} the buses to have PMU measurements, and moreover low latency communication links that convey the PMU measurements from all the buses to the control center that performs timely outage detection.

In practice, however, it may not be economically desirable that every bus has a PMU and constantly communicates the PMU measurements to the control center, \emph{particularly when there is very little performance loss with only a fraction of the buses providing PMU measurements}. We motivate the selection of a \emph{subset} of the buses to provide PMU measurements in the following two application scenarios:
\begin{enumerate}
\item As the cost of PMU installation in high voltage transmission networks is relatively high, we may want to install PMUs only at a subset of the buses to reduce cost.
\item Suppose over time the cost of PMUs drops and all (or many) buses have PMUs installed. Due to the \emph{information redundancy} in the PMU measurements, it may not be effective to communicate the data from all the buses to the control center, consuming an unnecessary amount of the expensive low-latency communication link capacity, and also incurring a longer data collection delay. Thus, the control center may select only a subset of the PMUs to provide measurements.
\end{enumerate}
In the following sections, we address the second application scenario, in which the PMU location selection can adapt to the changes of the power injections $\bm{P}$ and the network $\bm{B}$ matrix over time. We show later in section \ref{discgen} that the results can be directly extended to the first application scenario where a non-adaptive PMU location selection is required.
 
The problem of PMU location selection entails the following two questions:
\emph{
\begin{itemize}
\item Given that we want to choose $M~(2\le M\le N)$ of the $N$ buses to provide PMU measurements, which $M$ buses should we choose, and how can we characterize the corresponding outage detection performance? 
\item What is the optimal tradeoff between the number of PMUs $M$ and the outage detection performance?
\end{itemize}}
 
As mentioned in the last section, it makes no sense to have a bus chosen as the reference bus without having a PMU measuring its phase angle. Therefore, to find the optimal PMU locations, we cannot simply choose a reference bus arbitrarily, as it does lose generality. Instead, we \emph{traverse all the $N$ cases of choosing each bus as the reference bus.}

From now on, we assume that \emph{we have chosen a reference bus $r$ during this outer traversal, and optimize the location selection of the other $M-1$ PMUs.} We assume that all the phase angle vectors have been adjusted so that their $r^{th}$ entries always equal $0$:
\begin{align}
\theta^{(k)}_r = 0, \forall k = 0,1,\ldots,K. \nn
\end{align}
Accordingly, we denote by $\mc{M}$ the subset of buses with PMU measurements in addition to $r$, $\mc{M}\subseteq\{1,2,\ldots,N\}, |\mc{M}| = M-1$. As a result, for any two phase angle vectors, we cannot distinguish their entries at the other $N-M+1$ buses in $\mc{N}\backslash\mc{M}$. Therefore, for any phase angle vector $\bm{\theta}$, we project it into a sub-vector by extracting the $M-1$ entries of $\bm{\theta}$ whose indices are in $\mc{M}$ (it is not necessarily an orthogonal projection) and denote the projected vector by $\bm{\theta}_{\mc{M}} \in\mbb{R}^{(M-1)\times1}$. Accordingly, every set $\mc{M}$ leads to a set of projected signatures $\{\bm{\theta}_{\mc{M}}^{(k)}, k=0,1,\ldots,K\}$.

Clearly, with different choices of $\mc{M}$, the \emph{degrees of separation} among the projected signatures can differ considerably. Intuitively we want to optimize $\mc{M}$ to get better separation among the projected signatures, and thus better detection performance.

We model the separation among the projected signatures (and hence the outage detection performance) by 
the \emph{minimum distance} in $p$-norm among them:
\begin{align}
d_{\min}(\mc{M}) = \min_{0\le i< j\le K}\Vert\bm{\theta}_\mc{M}^{(i)} - \bm{\theta}_\mc{M}^{(j)}\Vert_p, \label{dmindef}
\end{align}
where $p$ is a parameter to choose, and $p=2$ corresponds to the Euclidian distance. Given $M$ as the total number of PMUs to use, 
$\mc{M}$ can then be optimized under the following \emph{Max-Min Distance Criterion}: 
\begin{align}
\max_{\mc{M},|\mc{M}| = M-1, r\notin\mc{M}}d_{\min}(\mc{M}). \label{maxminorig}
\end{align} 
Clearly, \eqref{maxminorig} is an NP hard combinatorial optimization of the set $\mc{M}$.

\subsection{An Integer Programming Formulation}
Note that $d_{\min}(\mc{M})$ can be re-written as
\begin{align}
d_{\min}(\mc{M}) = \left(\min(\bm{w}^T\bm{\Theta})\right)^{\frac{1}{p}},
\end{align}
where
\begin{itemize}
\item ~ \\
\vspace{-25pt} \begin{equation} \label{wdef} \bm{w}\in\mbb{R}^{N\times1},\text{ and }w_i = \left\{
\begin{array}{lll}
1, &\text{if } i\in\mc{M} \text{ or } i = r, \\
0, &\text{otherwise}.
\end{array}
\right.\end{equation}
\item $\bm{\Theta} \in\mbb{R}^{N\times {K+1 \choose 2}}$, and its columns are constructed by collecting the following ${K+1 \choose 2}$ $N\times1$ vectors:\\ $\forall 0\le i< j\le K,$ 
    \begin{align}
    \vert\bm{\theta}^{(i)} - \bm{\theta}^{(j)}\vert^p, \label{bigth}
    \end{align}
\end{itemize}
where the $\vert\cdot\vert^p$ operation is applied \emph{elementwise}.

In other words, the indices of the non-zero entries of $\bm{w}$ denote the $M$ buses that are chosen to provide PMU measurements (including the reference bus $r$). We name the binary vector $\bm{w}$ the \emph{bus selection indicator vector}.
Consequently, \eqref{maxminorig} is equivalent to the following integer programming problem:
\begin{align}
\max_{\bm{w}}&\min(\bm{w}^T\bm{\Theta}) \label{inlp}\\
s.t.~&w_i \in\{0,1\}, i=1,\ldots,N, \label{intc}\\
&\sum_{i=1}^{N}w_i = M, ~w_r = 1. \nn 
\end{align} 
Clearly, finding the global optimum of \eqref{inlp} requires a worst-case computational complexity of ${{N-1} \choose {M-1}}$.

\section{PMU Location Selection with Max-Min Distance Criterion}
In this section, we provide two algorithms for solving the combinatorial optimization of PMU location selection \eqref{maxminorig}. 
\subsection{A Greedy Heuristic}
We develop a greedy heuristic that generates a series of PMU location selection solutions $\mc{M}_2, \mc{M}_3,\ldots, \mc{M}_N$ for $M=2,3,\ldots,N$ respectively, that satisfy the following \emph{consistency property}:
\begin{align}
\mc{M}_2 \subset \mc{M}_3 \subset \ldots \subset \mc{M}_N. \label{conssol}
\end{align}
We present the greedy algorithm as follows in a slightly more general form with an arbitrary initial set of selected buses $\mc{M}_{ini}$:
\vspace{-10pt}

\begin{tabular}[c]{@{} l @{} l @{} p{8.4cm} @{}}
\\
\multicolumn{3}{@{}c}{\emph{Algorithm 1}: Greedy PMU Location Selection}\\
\vspace{-7pt}\\
\hline
\multicolumn{3}{@{}l}{Given $M$ and an initial set of buses $\mc{M}_{ini}$:}\\
\multicolumn{3}{@{}l}{~~~$m=|\mc{M}_{ini}|, \mc{M}_m = \mc{M}_{ini}$.}\\

\multicolumn{3}{@{}l}{Repeat}\\

&& $m\leftarrow m+1$,\\
\vspace{-25pt}\\
&& \begin{equation}\mc{M}_m = \mc{M}_{m-1} \cup \left\{\argmax_{n\in\mc{N}\backslash\mc{M}_{m-1}} d_{\min}\left(\mc{M}_{m-1}\cup \{n\}\right)\right\}\end{equation}\\
\vspace{-30pt}\\
\multicolumn{3}{@{}l}{Until $m=M$.} \\
\hline
\\
\end{tabular}
In other words, given the total number of PMUs $M$ and a starting set of chosen buses $\mc{M}_{ini}$, we choose another $M-|\mc{M}_{ini}|$ buses \emph{one by one}: At each step, we keep the already chosen buses; from the remaining buses, we choose the one that \emph{maximizes the current step's minimum distance}, and include it in the set of the chosen buses. 
When the only prior knowledge of the bus selection is the chosen reference bus $r$, $\mc{M}_{ini} = \{r\}$, and Algorithm 1 generates the set of greedy consistent solutions \eqref{conssol} as $M$ increases from $2$ to $N$.

\subsection{A Branch and Bound Algorithm}
First, we note that the integer programming formulation \eqref{inlp} has a concave objective function, and hence has the following \emph{relaxation} as a convex optimization:
\begin{align}
\max_{\bm{w}}&\min(\bm{w}^T\bm{\Theta}) \label{nlp}\\
s.t.~&0\le w_i\le1, i=1,\ldots,N, \label{linc}\\
&\sum_{i=1}^{N}w_i = M, ~w_r = 1. \nn
\end{align}
In fact, this convex optimization can be equivalently cast as a linear program \cite{BV04}. Accordingly, the optimal value of \eqref{nlp} serves as an \emph{upper bound}, denoted by $U_1$, on the global optimum of \eqref{inlp}. Meanwhile, Algorithm 1 with $\mc{M}_{ini} = \{r\}$ provides a \emph{lower bound}, denoted by $L_1$. 

\begin{RK}[A note on the rounding heuristic]
Another lower bounding heuristic is to find an integral solution by \emph{rounding} the fractional solution obtained from the relaxed problem \eqref{nlp}. In particular, we consider the heuristic by rounding \emph{the $M$ largest} fractional entries to $1$, and the others to $0$. In all scenarios that we simulated, we compared this rounding heuristic with the greedy heuristic: somewhat to our surprise, \emph{the greedy heuristic uniformly outperforms the rounding heuristic}. The reason is that, in the relaxed fractional solution, it is often some \emph{very small} non-zero entries that are \emph{critical} in the sense that losing them will drastically reduce the minimum distance. Consequently, for lower bounding the global optimum, we propose the greedy heuristic instead of the rounding heuristic, as the former is both much cheaper computationally and much better in performance. Thus, we use the relaxation technique \emph{not} to provide a fractional solution to round, but to \emph{upper bound} the global optimum and to develop a \emph{branch and bound} algorithm that can significantly improve the greedy solutions in a few iterations.
\end{RK}

For any bus $n\ne r$, \eqref{inlp} can be split into two sub-problems by fixing $w_n$ to be either $0$ or $1$:
\begin{align}
\max_{\bm{w}}&\min(\bm{w}^T\bm{\Theta}) \label{sub1}\\
s.t.~&w_i \in\{0,1\}, i=1,\ldots,N, \label{intc1}\\
&\sum_{i=1}^{N}w_i = M, ~w_n = 0, w_r = 1. \nn
\end{align}
\vspace{-10pt}
and
\vspace{-6pt}
\begin{align}
\max_{\bm{w}}&\min(\bm{w}^T\bm{\Theta}) \label{sub2}\\
s.t.~&w_i \in\{0,1\}, i=1,\ldots,N, \label{intc2}\\
&\sum_{i=1}^{N}w_i = M, ~w_n = 1, w_r = 1. \nn
\end{align}
Similarly to \eqref{nlp}, relaxations of these two sub-problems can be formed by replacing \eqref{intc1} and \eqref{intc2} with \eqref{linc}, and they provide upper bounds, denoted by $u_2^{(0)}$ and $u_2^{(1)}$, on the global optimum of \eqref{sub1} and \eqref{sub2} respectively. Meanwhile, applying the greedy heuristic under the constraint $w_n=0$ or $w_n=1$ provides lower bounds, denoted by $l_2^{(0)}$ and $l_2^{(1)}$, on these sub-problems' global optima. Define
\begin{align}
U_2\triangleq\max\{u_2^{(0)}, u_2^{(1)}\}, \text{ and } L_2\triangleq\max\{l_2^{(0)}, l_2^{(1)}\}. \label{updbound0}
\end{align}
Then, $U_2$ and $L_2$ are new upper and lower bounds on the original global optimum \eqref{inlp} \cite{BBB91}.

More generally, the above splitting procedure with relaxations and greedy heuristics can be applied on the sub-problems themselves to form more children sub-problems with upper and lower bounds. For example, for any bus $s, (s\ne n,r,)$ \eqref{sub1} can be further split into two sub-problems by adding yet another constraint $w_s=0$ or $w_s=1$ respectively.

We define the following upper and lower bounding oracles, as well as an oracle that returns the next bus to split:
\begin{define}
Oracle $U\!B(\mc{C})$ takes a \emph{constraint set} $\mc{C}$ as input, where $\mc{C}$ specifies a set of buses whose selection indicator variables are pre-determined to be either $0$ or $1$. An IP under the constraints $\mc{C}$ is formed, a relaxation is solved, and the optimum of this relaxation is output by $U\!B(\mc{C})$ as an upper bound on the optimum of the constrained IP.
\end{define}

For example, in \eqref{sub1} and \eqref{sub2}, the constraint sets are $\mc{C}^{(0)} = \{w_n = 0, w_r = 1\}$ and $\mc{C}^{(1)} = \{w_n = 1, w_r = 1\}$, respectively.

\begin{define}
Oracle $L\!B(\mc{C})$ takes a constraint set $\mc{C}$ as input. An IP under the constraints $\mc{C}$ is formed, a greedy solution is found by Algorithm 1, and the achieved objective value is output by $L\!B(\mc{C})$ as a lower bound on the optimum of the constrained IP.
\end{define}
\begin{define}
Based on the \emph{order} of the buses chosen by Algorithm 1, Oracle $next(\mc{C})$ outputs \emph{the first} bus that is chosen by this heuristic.
\end{define}

When a sub-problem with constraints $\mc{C}$ needs to be split further, $next(\mc{C})$ is the bus we choose to perform the splitting by fixing $w_{next(\mc{C})}$ to be either $0$ or $1$.

We now provide a branch and bound algorithm as in Algorithm 2 where $i_{\max}$ is the maximum number of iterations allowed. As the algorithm progresses, \emph{a binary tree is developed where each node represents a constraint set}. The \emph{leaf nodes} are kept in $\mc{S}$. The tree starts with a single node, $\{w_r = 1\}$, corresponding to the prior knowledge that the reference bus is $r$ who uses a PMU. When a sub-problem corresponding to a leaf node $\mc{C}^*$ is split into two new sub-problems, the two new constraint sets $\mc{C}^{(0)}$ and $\mc{C}^{(1)}$ become the children of the parent constraint set $\mc{C}^*$. 

In Algorithm 2, 
\eqref{updbound} is a generalization of \eqref{updbound0}. It means that the current global upper bound equals the \emph{highest} upper bound among all the \emph{leaf node} constraint sets. This is true because all the leaf nodes $\mc{S}$ represent a \emph{complete partition} of the original parameter space \cite{BBB91}. 
At the beginning of every iteration, in choosing which leaf node to split \eqref{maxubchoose}, we select the one that gives the \emph{highest} upper bound (i.e., the current \emph{global} upper bound). It is a heuristic based on the reasoning that, by further splitting this critical leaf node, a lower global upper bound may be obtained, (whereas splitting any other node will leave the global upper bound unchanged.) 
At iteration $i$, the current upper and lower bounds on the global optimum are available as $U_i$ and $L_i$. When these two bounds meet, i.e., $U_i - L_i < \epsilon$, the solution that achieves the current lower bound is guaranteed to be globally optimal. 

\begin{tabular}[c]{@{} l @{} l @{~} p{8.0cm} @{}}
\\
\multicolumn{3}{@{}c}{\emph{Algorithm 2:}} \\
\multicolumn{3}{@{}c}{PMU Location Selection using Branch and Bound}\\
\vspace{-7pt}\\
\hline
\multicolumn{3}{@{}l}{Initial step: $i=1$,} \\
\multicolumn{3}{@{}l}{~~~the initial constraint set: $~ \mc{C}_1 = \{w_r=1\}$,} \\
\multicolumn{3}{@{}l}{~~~the initial set of leaves of the tree of constraint sets}\\
\multicolumn{3}{@{}l}{~~~(initially a single node): $\mc{S} = \{\mc{C}_1\}$.} \\
&\multicolumn{2}{@{}l}{~~~Compute $U_1 = U\!B(\mc{C}_1),~L_1 = L\!B(\mc{C}_1)$.} \\

\multicolumn{3}{@{}l}{While $U_i-L_i>\epsilon$ or $i < i_{\max}$, repeat}\\

&~~~& Choose which leaf node constraint set to split: \\
\vspace{-20pt}\\
&~~~& \begin{equation}\mc{C}^* = \argmax_{\mc{C}\in\mc{S}}\{U\!B(\mc{C})\}. \label{maxubchoose}\end{equation}\\
\vspace{-20pt}\\
&~~~& Choose the next bus to split, $n = next(\mc{C}^*)$,\\
&~~~& Form two new constraint sets, \\
&~~~& ~~~$\mc{C}_{i+1}^{(0)} = \mc{C^*}\cup\{w_n=0\},~ \mc{C}_{i+1}^{(1)} = \mc{C^*}\cup\{w_n=1\}$.\\
&~~~& In the set of leaves $\mc{S}$, replace the parent constraint set $\mc{C}^*$ with the two children $\mc{C}_{i+1}^{(0)}$ and $\mc{C}_{i+1}^{(1)}$:\\
&~~~& ~~~$\mc{S} \leftarrow \left(\mc{S}\backslash\{\mc{C}^*\}\right)\cup \{\mc{C}_{i+1}^{(0)}\}\cup \{\mc{C}_{i+1}^{(1)}\}$.\\

&~~~& Compute new upper and lower bounds for the two new constrained IP: \\
&~~~& ~~~$U\!B(\mc{C}_{i+1}^{(0)}),~ U\!B(\mc{C}_{i+1}^{(1)}),~ L\!B(\mc{C}_{i+1}^{(0)}),~ L\!B(\mc{C}_{i+1}^{(1)})$, \\
&~~~& Update the global upper and lower bounds, \\
\vspace{-20pt}\\
&~~~& \begin{equation} U_{i+1} = \max_{\mc{C}\in\mc{S}}\{U\!B(\mc{C})\}, ~L_{i+1} = \max_{\mc{C}\in\mc{S}}\{L\!B(\mc{C})\}. \label{updbound}\end{equation} \\
\vspace{-20pt}\\
\multicolumn{3}{@{}l}{$i\leftarrow i+1$.} \\
\vspace{-5pt}\\
\multicolumn{3}{@{}l}{Choose the best achieved solution so far:} \\
\multicolumn{3}{@{}l}{~~~~$\hat{\mc{C}} = \argmax_{\mc{C}\in\ms{S}}L\!B(\mc{C})$.} \\
\multicolumn{3}{@{}l}{~~~~Return the greedy solution under the constraint set $\hat{\mc{C}}$.} \\
\hline
\\
\end{tabular}

We note that Algorithm 1 is a degraded version of Algorithm 2 with just one iteration. As the total number of possible constraint sets is $2^N$ (corresponding to the $2^N$ bus selection indicator vectors $\bm{w}$), Algorithm 2 is guaranteed to converge in $2^N$ iterations (and in practice much less as will be shown later.) To limit the algorithm's run time, a maximum number of iterations $i_{\max}$ can be enforced as in Algorithm 2.

Finally, we define $i_\text{achieve}$ to be the number of iterations used to \emph{achieve} the globally optimal solution, and $i_\text{prove}$ the number of iterations used to \emph{prove} its global optimality. In other words, it takes $i_\text{achieve}$ iterations for the \emph{lower bound} to reach the global optimum, while it takes $i_\text{prove}$ iterations for \emph{both} the upper and lower bounds to reach the global optimum. As will be shown next, typically we have $i_\text{achieve}\ll i_\text{prove}$.

\section{Performance Evaluation and the Optimal $M$-$d_{\min}$ Tradeoff}
In this section, we simulate the proposed algorithms in IEEE $14, 24$, and $30$ bus systems using the software toolbox MATPOWER \cite{MatP11}. We set the outage event set $\mc{E}$ to be the set of all \emph{single line outages} that do not create islands. 
Using Algorithm 2, we find the \emph{globally optimal} solutions for PMU location selection, and characterize the \emph{optimal tradeoff} between the number of PMUs used and the minimum distance in Euclidean norm among the projected outage signatures.

\begin{figure}[tb]
  \centering
  \includegraphics[scale=0.58]{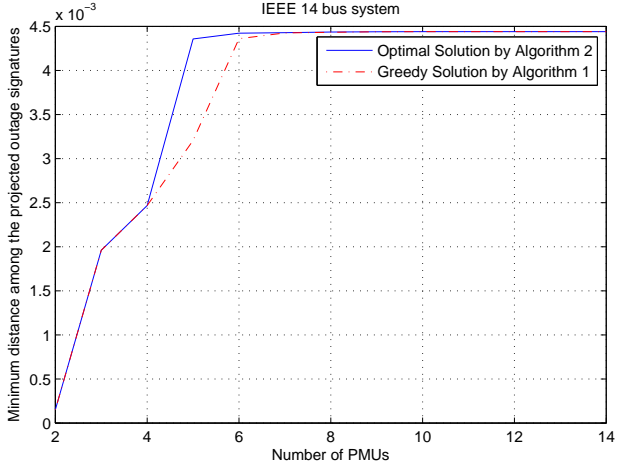}
  \caption{Tradeoff between the number of PMUs and the maximum achievable $d_{\min}$ among the projected outage signatures; IEEE 14 bus system.}
  \label{dmin14opt}
\end{figure}

\begin{figure}[tb]
  \centering
  \includegraphics[scale=0.58]{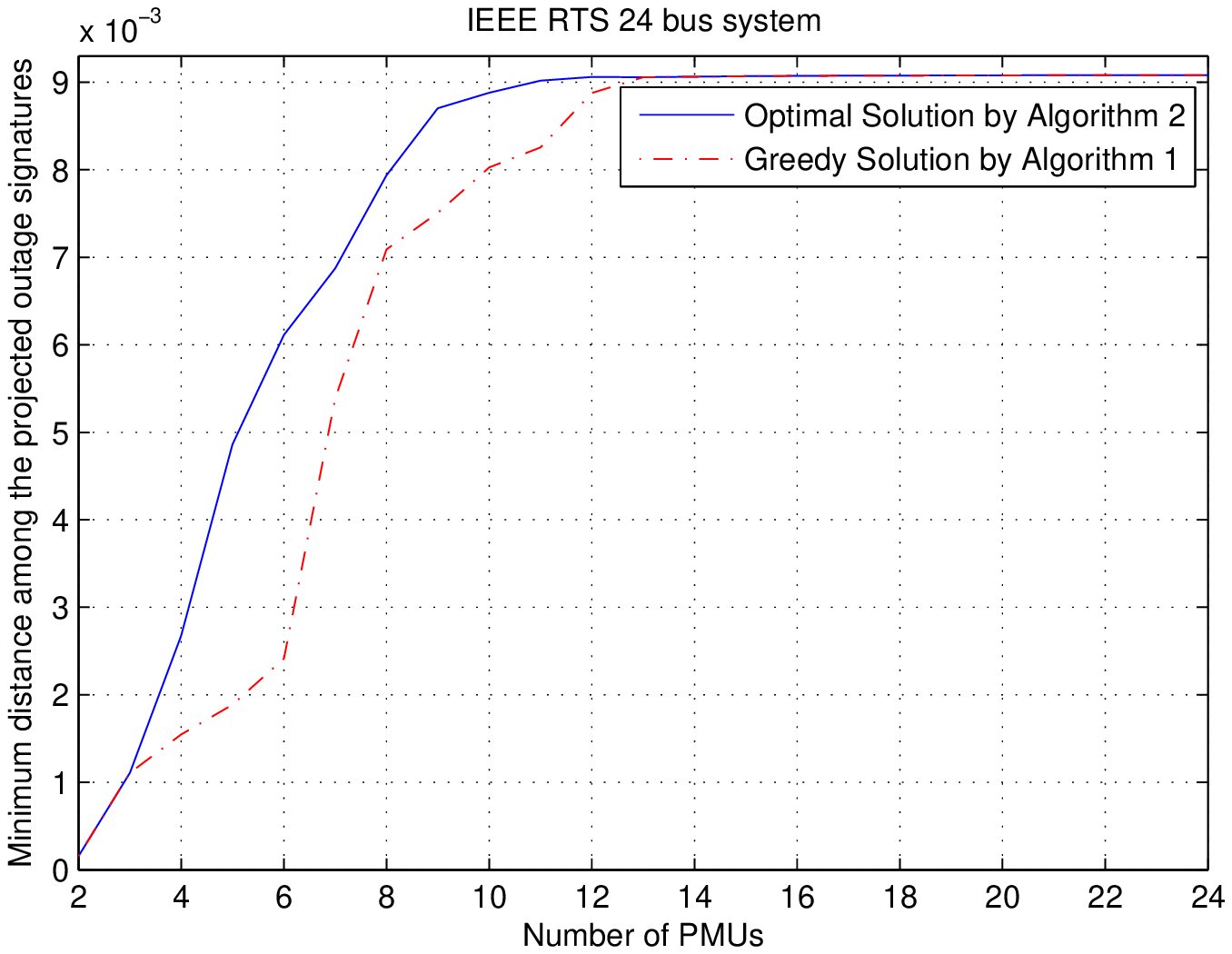}
  \caption{Tradeoff between the number of PMUs and the maximum achievable $d_{\min}$ among the projected outage signatures; IEEE RTS 24 bus system.}
  \label{dmin24opt}
\end{figure}

\begin{figure}[tb]
  \centering
  \includegraphics[scale=0.58]{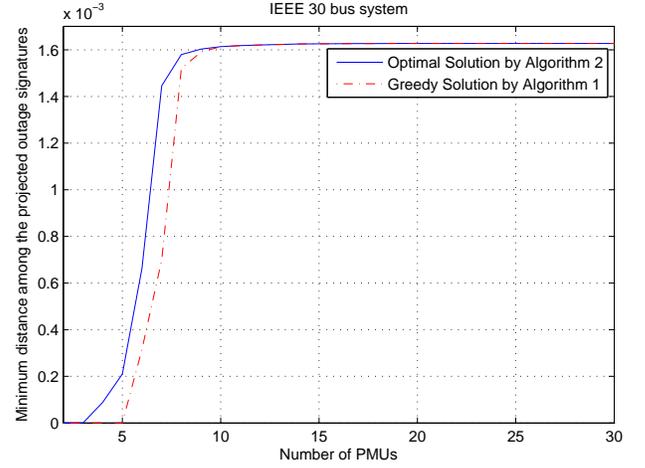}
  \caption{Tradeoff between the number of PMUs and the maximum achievable $d_{\min}$ among the projected outage signatures; IEEE 30 bus system.}
  \label{dmin30opt}
\end{figure}

We denote the maximum achievable minimum distance as a function of $M$ by
\begin{align}
d_{\min}^*(M), M=2,3,\ldots,N.
\end{align}
We plot $d_{\min}^*(M)$ in Figure \ref{dmin14opt}, \ref{dmin24opt} and \ref{dmin30opt}, for the 14, 24 and 30 bus systems respectively. In comparison, the greedy solutions achieved by Algorithm 1 are also plotted.

We make the following observations:
\begin{itemize}
\item Having $M \approx \frac{N}{3}$ PMUs is roughly sufficient to achieve the same optimal $d_{\min}$ as with $M=N$ (i.e., all the buses having PMU measurements).
\item As $M$ increases from 2, every additional PMU significantly increases $d_{\min}^*(M)$, until roughly $d_{\min}^*(N)$ is achieved (with $M \approx \frac{N}{3}$ as mentioned above).
\item When the available number of PMUs is sufficiently large ($M$ greater than $\frac{N}{3}\sim\frac{N}{2}$ in these examples), globally optimal performance can be achieved by simple greedy solutions.
\end{itemize}

The number of iterations needed for Algorithm 2 to reach the globally optimal solutions, $i_\text{achieve}$ and $i_{\text{prove}}$, are summarized in Table 1, 2 and 3, for the 14, 24 and 30 bus systems respectively.

\begin{table}[!bp]
\begin{tabular}[l]{|p{0.8cm}|p{0.4cm}|p{0.4cm}|p{0.4cm}|p{0.4cm}|p{0.4cm}|p{0.4cm}|p{0.4cm}|p{0.4cm}|p{0.4cm}|} 
\multicolumn{10}{@{}c}{\emph{Table 1}: Number of iterations to reach the global optimum, 14 bus.}\\
\hline
$M$ & 2 & 3 & 4 & 5 & 6 & 7 & 8 & 9 & 10 \\
\hline
$i_\text{achieve}$ & 1 & 1 & 1 & 17 & 2 & 1 & 1 & 1 & 1\\
$i_{\text{prove}}$ & 3 & 5 & 16 & 17 & 2 & 1 & 1 & 1 & 1\\
\hline
$M$ & 11 & 12 & 13 & 14 & \multicolumn{5}{|c|}{ } \\ 
\hline
$i_\text{achieve}$ & 1 & 1 & 1 & 1 & \multicolumn{5}{|c|}{ }\\
$i_{\text{prove}}$ & 1 & 1 & 1 & 1 & \multicolumn{5}{|c|}{ }\\
\hline
\end{tabular}
\end{table}

\vspace{5pt}

\begin{table}[!bp]
\begin{tabular}[c]{|p{0.8cm}|p{0.4cm}|p{0.4cm}|p{0.4cm}|p{0.4cm}|p{0.4cm}|p{0.4cm}|p{0.4cm}|p{0.4cm}|p{0.4cm}|} 
\multicolumn{10}{@{}c}{\emph{Table 2}: Number of iterations to reach the global optimum, 24 bus.}\\
\hline
$M$ & 2 & 3 & 4 & 5 & 6 & 7 & 8 & 9 & 10 \\ 
\hline
$i_\text{achieve}$ & 1 & 1 & 3 & 3 & 3 & 17 & 6 & 6 & 12\\
$i_{\text{prove}}$ & 2 & 12 & 40 & 83 & 144 & 395 & 268 & 208 & 171\\
\hline
$M$ & 11 & 12 & 13 & 14 & 15 & 16 & 17 & 18 & 19 \\
\hline
$i_\text{achieve}$ & 5 & 3 & 1 & 1 & 1 & 1 & 1 & 1 & 1 \\
$i_{\text{prove}}$ & 170 & 3 & 1 & 1 & 1 & 1 & 1 & 1 & 1 \\
\hline
$M$ & 20 & 21 & 22 & 23 & 24 & \multicolumn{4}{|c|}{ } \\ 
\hline
$i_\text{achieve}$ & 1 & 1 & 1 & 1 & 1 & \multicolumn{4}{|c|}{ }\\
$i_{\text{prove}}$ & 1 & 1 & 1 & 1 & 1 & \multicolumn{4}{|c|}{ }\\
\hline
\end{tabular}
\end{table}

\vspace{5pt}

\begin{table}[!bp]
\begin{tabular}[c]{|p{0.8cm}|p{0.4cm}|p{0.4cm}|p{0.4cm}|p{0.4cm}|p{0.4cm}|p{0.4cm}|p{0.4cm}|p{0.4cm}|p{0.4cm}|p{0.4cm}|}
\multicolumn{10}{@{}c}{\emph{Table 3}: Number of iterations to reach the global optimum, 30 bus.}\\
\hline
$M$ & 2 & 3 & 4 & 5 & 6 & 7 & 8 & 9 & 10\\
\hline
$i_\text{achieve}$ & 1 & 1 & 19 & 4 & 5 & 5 & 4 & 5 & 1 \\
$i_{\text{prove}}$ & 3 & 21 & 99 & 53 & 24 & 30 & 31 & 5 & 1 \\
\hline
$M$ & 11 & 12 & 13 & 14 & 15 & 16 & 17 & 18 & 19\\
\hline
$i_\text{achieve}$ & 1 & 1 & 1 & 1 & 1 & 1 & 1 & 1 & 1 \\
$i_{\text{prove}}$ & 1 & 1 & 1 & 1 & 1 & 1 & 1 & 1 & 1 \\
\hline
$M$ & 20 & 21 & 22 & 23 & 24 & 25 & 26 & 27 & 28\\
\hline
$i_\text{achieve}$ & 1 & 1 & 1 & 1 & 1 & 1 & 1 & 1 & 1 \\
$i_{\text{prove}}$ & 1 & 1 & 1 & 1 & 1 & 1 & 1 & 1 & 1 \\
\hline
$M$ & 29 & 30 & \multicolumn{7}{|c|}{ }\\
\hline
$i_\text{achieve}$ & 1 & 1 & \multicolumn{7}{|c|}{ }\\
$i_{\text{prove}}$ & 1 & 1 & \multicolumn{7}{|c|}{ }\\
\hline
\end{tabular}
\end{table}

\begin{figure}[tb]
  \centering
  \includegraphics[scale=0.58]{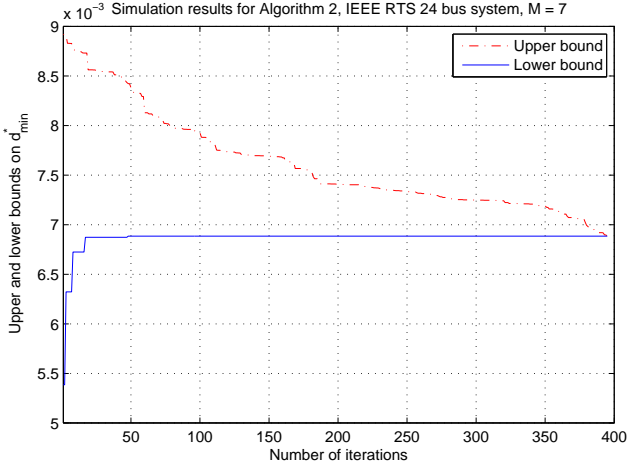}
  \caption{The instant upper and lower bounds on the maximum achievable minimum distance as Algorithm 2 iterates; IEEE RTS 24 bus system, $M=7$.}
  \label{24busM7}
\end{figure}

We make the following observations on the efficiency of Algorithm 2:
\begin{itemize}
\item In these three systems, Algorithm 2 always finds the globally optimal solution in \emph{less than 19 iterations}. For solving the $M=4$ case in the 30 bus system, $i_\text{achieve}=19$. For most other cases, $i_\text{achieve}$ is much smaller.
\item When $M$ is relatively small, it takes a much larger number of iterations to \emph{prove} that the solutions achieved within $19$ iterations are indeed globally optimal. The maximum $i_\text{prove}$ is 395 for solving the $M=7$ case in the 24 bus system. For this particular case, the upper and lower bounds achieved as Algorithm 2 progresses are plotted in Figure \ref{24busM7}.
\item $i_\text{achieve}$ and $i_\text{prove}$ are much smaller than the \emph{combinatorial complexity} of this NP hard problem. E.g., for solving the case of $M = 7$ in the 24 bus system, an exhaustive search has to traverse ${24 \choose 7} = 346104$ PMU location selections, whereas Algorithm 2 finds the optimal solution in $i_\text{achieve} = 17$ iterations, and verifies its optimality in $i_\text{prove} = 395$ iterations.
\end{itemize}
The fact that $i_\text{achieve}$ is typically small demonstrates that Algorithm 2 is very efficient in finding the globally optimal PMU locations. We note that, in practice, the sometimes larger $i_\text{prove}$ does not matter at all. This is because we always use the best solution found within a predetermined run time, and as long as $i_\text{achieve}$ is sufficiently small, we will find the optimal solution (albeit without proving its optimality).

\begin{figure}[tb]
  \centering
  \includegraphics[scale=0.35]{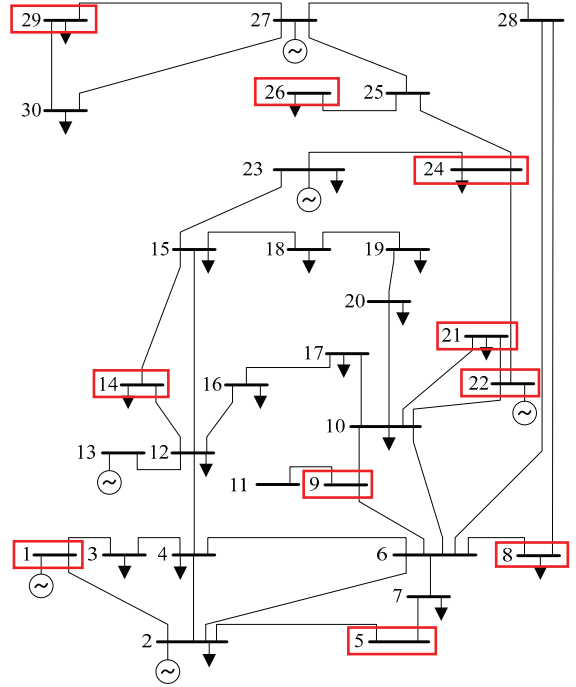}
  \caption{One line diagram of the IEEE 30 bus system, and the optimal PMU location selection with $M = 10$. The buses enclosed by rectangles are the optimally chosen buses.}
  \label{30busM10}
\end{figure}

Finally, as an example of the optimal PMU location selection, we depict the buses that are chosen for the 30 bus system with $M = 10$ (cf. Figure \ref{30busM10}). In this case, bus $1,5,8,9,14,21,22,24,26,29$ are chosen to provide PMU measurements, achieving the same minimum distance among the outage signatures as with all the 30 buses chosen. As shown in Table 3, since $i_\text{achieve} = 1$ for $M=10$, this globally optimal solution can in fact be found by the simple greedy method of Algorithm 1.

\section{Discussion}
\subsection{Normalizing Distances by Noise Variances} \label{gennoise}
The max-min distance criterion we have used is for countering the effect of the noise \eqref{noisyv} by separating the outage signatures as much as possible. In general, the variances of the noise in the PMU measurements \emph{may differ at different buses}. This shall be taken into account when defining the \emph{effective distances} among the outage signatures.

For two outage signatures $\bm{\theta}^{(i)}$ and $\bm{\theta}^{(j)}$, denote their difference in absolute value by $\Delta\bm{\theta} = |\bm{\theta}^{(i)} - \bm{\theta}^{(j)}|$, where $|\cdot|$ is applied elementwise. The $n^{th}$ entry $\Delta\theta_n$ is the distance between the two signatures \emph{at the $n^{th}$ bus}.

Suppose now that the noises are independently Gaussian with variances $\sigma_1^2, \sigma_2^2, \ldots, \sigma_N^2$ at the $N$ buses. We can then \emph{normalize} the distance at the $n^{th}$ bus by a factor of $\frac{1}{\sigma_n}$:
\begin{align}
\Delta\tilde{\bm{\theta}} = \bm{\Sigma}^{-1}\Delta\bm{\theta},
\end{align}
where $\bm{\Sigma} = \diag(\sigma_1, \sigma_2, \ldots,\sigma_N)$. After this normalization, it is now appropriate to \emph{add up} all the buses' contributions in separating the signatures by computing the $p$-norm of $\Delta\tilde{\bm{\theta}}$:
\begin{align}
\Vert\Delta\tilde{\bm{\theta}}\Vert_p = \left(\sum_{n=1}^N\Delta\tilde{\theta}_n^p\right)^\frac{1}{p}.
\end{align}

Accordingly, the definition of the minimum distance among the signatures \eqref{dmindef} is generalized as follows:
\begin{align}
\tilde{d}_{\min}(\mc{M}) 
= \min_{0\le i< j\le K}\left\Vert\left(\bm{\Sigma}^{-1}\left(\bm{\theta}^{(i)} - \bm{\theta}^{(j)}\right)\right)_\mc{M}\right\Vert_p \nn
\end{align}

For the IP formulation \eqref{inlp}, it was assumed that $\bm{\Sigma} = \sigma \bm{I}$ for some $\sigma$ in the previous definition of $\bm{\Theta}$ \eqref{bigth}. With an arbitrary $\bm{\Sigma}$, the matrix of the distances among the outage signatures is then generalized by the following normalization:
\begin{align}
\tilde{\bm{\Theta}} = \bm{\Sigma}^{-p}\bm{\Theta}.
\end{align}


\subsection{Extension to Non-Adaptive PMU Location Selection} \label{discgen}
We have provided two application scenarios in Section \ref{pmuloc} that motivate the problem of PMU location selection. 
In the previous sections, we have focused on the second scenario in which the outage signatures $\{\bm{\theta}^{(k)},k=0,1,\ldots,K\}$, and hence the optimal PMU location selection, depend on the current network setting $\bm{P}$ and $\bm{B}$.

In the first scenario, however, location selection for \emph{PMU installation} is a planning problem whose solution must work \emph{non-adaptively} in the sense that the installed PMUs are not transferable to other locations. As a result, the PMU locations must be chosen to accommodate, if not all, \emph{the most typical} network settings in terms of $\bm{P}$ and $\bm{B}$. As we consider wide-area transmission networks, this requirement is addressed as follows:
\begin{itemize}
    \item
     The network topology $\bm{B}$ is in general slowly and slightly changing, and an approximate estimate of it would suffice for computing typical sets of outage signatures.
    \item Typical power injections $\bm{P}$ are available from grid statistics for high voltage transmission networks: collecting the typical data for \emph{different seasons} and during \emph{daytime and nighttime} would suffice for computing typical sets of outage signatures.
\end{itemize}

With a set of typical power injections $\bm{P}_1, \bm{P}_2, \ldots, \bm{P}_T$, and an estimated $\bm{B}$, there are $T$ typical sets of outage signatures:
\begin{align}
\{\bm{\theta}^{(k,t)}, k = 0,1,\ldots,K\}, t = 1,2,\ldots,T.
\end{align}
The minimum distance \eqref{dmindef} can be generalized as
\begin{align}
d_{\min}(\mc{M}) = \min_{1\le t\le T}\min_{0\le i< j\le K}\Vert\bm{\theta}_\mc{M}^{(i,t)} - \bm{\theta}_\mc{M}^{(j,t)}\Vert_p. \label{dmindefgen}
\end{align}
Similarly to \eqref{bigth}, we construct $\bm{\Theta} \in\mbb{R}^{N\times T\cdot{K+1 \choose 2}}$, and its columns are constructed by collecting the following $T\cdot{K+1 \choose 2}$ $N\times1$ vectors: $\forall 1\le t\le T, ~\forall 0\le i< j\le K,$ 
    \begin{align}
    \vert\bm{\theta}^{(i,t)} - \bm{\theta}^{(j,t)}\vert^p, \label{bigthgen}
    \end{align}
where the $\vert\cdot\vert^p$ operation is applied \emph{elementwise}.

The problem of non-adaptive PMU location selection with \eqref{dmindefgen} can then be formulated as the same integer programming problem \eqref{inlp}, for which Algorithm 1 and 2 can be applied.

\section{Concluding Remarks}
We have studied the problem of PMU location selection for detecting line outages in wide-area transmission networks. We have first introduced the voltage phase angle signatures of the outages. We have then established the PMU location selection problem as maximizing the minimum distance among the set of outage signatures, and have shown that it can be equivalently formulated as an integer programming problem. For this NP hard problem, we have developed a greedy heuristic and a linear programming relaxation, providing a lower and an upper bound on the global optimum respectively. Based on this heuristic and relaxation, we have proposed a branch and bound algorithm to find the globally optimal PMU locations. Using this algorithm, we have characterized the optimal trade-offs between the number of PMUs and the minimum distance among the outage signatures in IEEE 14, 24 and 30 bus systems. For all the simulated cases, the optimal PMU locations are found in at most $19$ (and in most cases a much less number of) iterations. From the optimal trade-offs, we have observed that it is sufficient to have roughly one third of the buses to provide PMU measurements in order to achieve the same outage detection performance as with all the buses providing PMU measurements.

For future work, while we have employed the max-min distance criterion among the outage signatures, it is interesting to investigate the use of other criteria in characterizing outage detection performance. We would also like to examine the optimal PMU location selections with outage event sets more general than line outages. Finally, it remains an interesting open question how to take into account the needs from \emph{multiple} tasks simultaneously (e.g., outage detection \emph{and} state estimation) while selecting PMU locations.


\bibliographystyle{IEEEbib}
{\bibliography{PES12}}

%








\end{document}